\documentclass[11pt, psamsfonts]{amsart}
\usepackage{amssymb, amsmath, amsthm}
\usepackage[final, hypertex]{hyperref}
\usepackage[a4paper, centering]{geometry}
\newtheorem{thm}{Theorem}[section]
\newtheorem{prop}[thm]{Proposition}
\newtheorem{lem}[thm]{Lemma}

\theoremstyle{definition}
\newtheorem{defn}[thm]{Definition}
\theoremstyle{remark}
\newtheorem{rem}[thm]{Remark}

\def\R{\mathbb{R}}

\def\N{\mathbb{N}}
\def\U{\mathcal{U}}

\def\pscal#1#2{\left\langle#1,\,#2\right\rangle}

\def\convopen{\mathcal{K}^n_0}

\def\gau#1{\rho_{#1}}
\def\gauge{\rho}
\def\pgauge{\gauge^0}

\def\distb#1{d_{#1}}
\def\dist{\distb{\Omega}}

\def\cut{m}
\def\uv{\nu}
\def\curv{\kappa}
\def\curvg{\tilde{\kappa}}
\def\nor{\nu}

\def\Lip{\textrm{Lip}}
\def\bw{\overline{W}}
\def\bv{\overline{V}}
\def\vf{v_f}

\def\matr{\mathcal{M}}
\def\uu{w}
\def\vv{z}
\def\vz{v}
\def\uz{u}
\def\Cb{C_b}
\def\wbot{\overline{w}}

\DeclareMathOperator{\spt}{supp}

\DeclareMathOperator{\dive}{div}
\DeclareMathOperator{\trace}{Tr}

\DeclareMathOperator{\proj}{\Pi}

\begin{document}
\title[PDEs of Monge-Kantorovich type]%
{On a system of partial differential equations
of Monge-Kantorovich type}

\author[G.~Crasta]{Graziano Crasta}
\address{Dipartimento di Matematica ``G.\ Castelnuovo'', Univ.\ di Roma I\\
P.le A.\ Moro 2 -- 00185 Roma (Italy)}
\email[Graziano Crasta]{crasta@mat.uniroma1.it}

\author[A.~Malusa]{Annalisa Malusa}
\email[Annalisa Malusa]{malusa@mat.uniroma1.it}



\begin{abstract}
We consider a system of PDEs of Monge-Kantorovich type
arising from
models in granular matter theory
and in electrodynamics of hard superconductors.
The existence of a solution of such system
(in a regular open domain $\Omega\subset\mathbb{R}^n$),
whose construction is
based on an asymmetric Minkowski distance from the
boundary of $\Omega$, was already established
in
[{G.} Crasta and {A.} Malusa, \emph{The distance function from the boundary in a
  {M}inkowski space}, to appear in Trans.\ Amer.\ Math.\ Soc.].
In this paper we prove that this solution is
essentially unique.
A fundamental tool in our analysis is
a new regularity result for an elliptic nonlinear equation
in divergence form,
which is of some interest by itself.
\end{abstract}

\maketitle

\section{Introduction}

Let $\gauge\colon\R^n\to\R$ be a $C^2$ gauge function,
i.e.\ a convex and positively $1$-homogenous function,
of class $C^2$ in $\R^n\setminus\{0\}$.
In this paper we are concerned with
the system of partial differential equations
\begin{equation}\label{f:syst1}
\begin{cases}
-\dive(v\, D\gauge(Du)) = f
&\textrm{in $\Omega$},\\
\gauge(Du)\leq 1
&\textrm{in $\Omega$},\\
\gauge(Du) = 1
&\textrm{in $\{v>0\}$},
\end{cases}
\end{equation}
complemented with the conditions
\begin{equation}\label{f:syst2}
\begin{cases}
u\geq 0,\
v\geq 0
&\textrm{in $\Omega$},\\
u=0
&\textrm{on $\partial\Omega$}.
\end{cases}
\end{equation}
Here $\Omega\subset\R^n$ is a bounded domain of class $C^2$
and $f\geq 0$ is a bounded continuous function in $\Omega$.
A solution of this system is a pair $(u,v)$ of nonnegative functions,
with $u$ Lipschitz continuous in $\overline{\Omega}$ and
$v$ bounded and continuous in $\Omega$, satisfying the following
additional conditions:
(a) $u = 0$ on $\partial\Omega$;
(b) $\gauge(Du)\leq 1$ almost everywhere in $\Omega$;
(c) $u$ is a viscosity solution of
$\gauge(Du) = 1$ in the open set $\{v>0\}$;
(d) $v$ is a solution
of the first equation in (\ref{f:syst1})
in the sense of distributions
(see Definition~\ref{d:sol} below).

This system of PDEs arises in some different situations.
For example,
the functions $u$ and $v$ can be interpreted respectively as
the magnetic field and the power dissipation
in a cylindrical hard superconductor of cross-section $\Omega$
exposed to an external magnetic field linearly increasing in time
(see e.g.~\cite{Cha}).
Moreover, in the case $\gauge(\xi) = |\xi|$, (\ref{f:syst1})-(\ref{f:syst2})
gives the stationary solutions of models in granular matter
theory (see \cite{CaCa,CCCG}).
Another application concerns the existence of solutions to nonconvex
minimum problems in calculus of variations
(see \cite{CCCG,CPT,Ce1}).
Finally, Bouchitt\'e and Buttazzo \cite{BoBu}
have studied a more general system
in order to describe optimal solutions
of some shape optimization problems.

The results presented in this paper are an extension
of those proved in \cite{CCCG}, where the
case $\gauge(\xi) = |\xi|$ was considered.
An explicit solution to (\ref{f:syst1})-(\ref{f:syst2}) was constructed
in \cite{CMf}.
In order to describe this solution, we need some additional notation.
Let $\dist\colon\overline{\Omega}\to\R$ be the Minkowski distance from
the boundary $\partial\Omega$, defined by
\[
\dist(x) = \inf_{y\in\partial\Omega} \pgauge(x-y),
\qquad x\in\overline{\Omega},
\]
where $\pgauge$ is the polar function of $\gauge$.
It is well-known that $\dist$ is Lipschitz continuous in
$\overline{\Omega}$, and it is the unique viscosity solution
of $\gauge(D\dist) = 1$ in $\Omega$ vanishing on $\partial\Omega$
(see \cite{Li}).
In \cite{CMf} it was shown that there exists
a bounded continuous function $\vf\colon\Omega\to [0,+\infty)$,
whose explicit expression depends on $f$, $\gauge$ and on the geometry
of $\Omega$ (see Section~\ref{s:exi}),
such that the pair $(\dist, \vf)$ is a solution to
(\ref{f:syst1})-(\ref{f:syst2}).

The aim of this paper is to show that this solution is essentially unique.
More precisely, we shall show that
if $(u,v)$ is a solution to (\ref{f:syst1})-(\ref{f:syst2}),
then $v=\vf$ and $u=\dist$ in
$\Omega_f = \{x\in\Omega;\ \vf(x) > 0\}$.
The proof of this uniqueness result is based on
several ingredients.
Some of them are an adaptation to our setting
of arguments developed in
\cite{Egan,CaCa,Pr,CCCG}.
A key point of the uniqueness proof consists in showing
that, if $(\dist, v)$ is a solution to
(\ref{f:syst1})-(\ref{f:syst2}), then $v$ vanishes
on the singular set $\Sigma$ of $\dist$
(see Proposition~\ref{vzero} below).
In this respect, we use here a blow-up argument
introduced by Evans and Gangbo
for the case $\gauge(\xi) = |\xi|$
(see \cite{Egan}, Section~7, and \cite{CCCG}),
which in turn relies on the regularity of the solutions
to the classical Laplace equation
$\Delta u = 0$ in an open set $A\subset\R^n$.
In our setting the classical Laplace equation is
replaced by
\begin{equation}\label{f:equa2b}
-\dive(D\gauge(Du)) = 0
\qquad \textrm{in}\ A\,.
\end{equation}
Since the function $a(\xi) := D\gauge(\xi)$ is defined
and positively $0$-homogeneous in $\R^n\setminus\{0\}$,
no standard regularity result can be applied.
In Section~\ref{s:reg} we prove that,
if $\gauge\in C^2(\R^n\setminus\{0\})$ and
$u$ is a Lipschitz continuous solution of (\ref{f:equa2b})
satisfying
$\gauge(Du) = 1$ almost everywhere in $A$,
then $u$ is of class $C^{1,\alpha}$ locally in $A$
(see Theorem~\ref{t:reg}).
Thanks to this regularity result,
the blow-up argument of Evans and Gangbo still
works in our setting
(see Proposition~\ref{vzero}).

\section{Notation and Preliminaries}

\subsection{Basic notation}
The standard scalar product of two
vectors $x,y\in\R^n$
is denoted by $\pscal{x}{y}$,
and $|x|$ denotes the
Euclidean norm of $x\in\R^n$.
Given two vectors $v,w\in\R^n$, the symbol
$v\otimes w$ will denote their tensor product,
i.e.\  the linear application from
$\R^n$ to $\R^n$ defined by
$(v\otimes w)(x) = v\,\pscal{w}{x}$.

By $S^{n-1}$ we denote the set of unit vectors
of $\R^n$, and by $\matr_k$ the set of
$k\times k$ square matrices.
We shall denote by $(e_1,\ldots,e_n)$
the standard basis of $\R^n$.
The closed segment joining $x\in\R^n$ to $y\in\R^n$
will be denoted by $[x,y]$,
while $(x,y)$ will denote the same segment without
the endpoints.

As customary, $B_r(x_0)$ and $\overline{B}_r(x_0)$
are respectively the open and the closed ball
centered at $x_0$ and with radius $r>0$.

Given $A\subset\R^n$,
we shall denote by $\Lip(A)$, $C(A)$, $\Cb(A)$ and $C^k(A)$, $k\in\N$
the set of functions $u\colon A\to\R$
that are respectively Lipschitz continuous, continuous,
bounded and continuous,
and $k$-times continuously differentiable in $A$.
(Here and thereafter $\N$ will denote the set of nonnegative
integers.)
Moreover, $C^{\infty}(A)$ will denote the set of functions
of class $C^k(A)$ for every $k\in\N$,
while $C^{k,\alpha}(A)$ will be the set of functions
of class $C^k(A)$ with H\"older continuous $k$-th partial derivatives
with exponent $\alpha\in [0,1]$.

A bounded open set $A\subset\R^n$
(or, equivalently, its closure
$\overline{A}$ or its boundary $\partial A$)
is of class $C^k$, $k\in\N$,
if for every point $x_0\in\partial A$
there exists a ball $B=B_r(x_0)$ and a one-to-one
mapping $\psi\colon B\to D$ such that
$\psi\in C^k(B)$, $\psi^{-1}\in C^k(D)$,
$\psi(B\cap A)\subseteq\{x\in\R^n;\ x_n > 0\}$,
$\psi(B\cap\partial A)\subseteq\{x\in\R^n;\ x_n = 0\}$.
If the maps $\psi$ and $\psi^{-1}$ are of class
$C^{\infty}$ or $C^{k,\alpha}$ ($k\in\N$, $\alpha\in [0,1]$),
then $A$ is said to be of class
$C^{\infty}$ or $C^{k,\alpha}$ respectively.

\subsection{Convex geometry}
By $\convopen$ we denote the class of
nonempty, compact, convex subsets of $\R^n$
with the origin as an interior point.
We shall briefly refer to the elements of $\convopen$
as \textsl{convex bodies}.
The polar body of a convex body $K\in\convopen$
is defined by
\[
K^0 = \{p\in\R^n;\ \pscal{p}{x}\leq 1\ \forall x\in K\}\,.
\]
We recall that, if $K\in\convopen$, then
$K^0\in\convopen$ and $K^{00} = (K^0)^0 = K$
(see \cite[Thm.~1.6.1]{Sch}).

Given $K\in\convopen$ we define
its gauge function as
\[
\gau{K}(\xi) = \inf\{ t\geq 0;\ \xi\in t K\}\,.
\]
It is easily seen that
\[
\gau{K^0}(\xi) =
\sup\left\{\pscal{\xi}{p};\ p\in K\right\}\,,
\]
i.e.\  the gauge function of the polar set $K^0$
coincides with the support function of the set $K$.
Let $0<c_1\leq c_2$ be such that
$\overline{B}_{c_2^{-1}}(0)\subseteq K \subseteq \overline{B}_{c_1^{-1}}(0)$.
Upon observing that $\xi/\gau{K}(\xi)\in K$ for every $\xi\neq 0$,
we get
\begin{equation}\label{f:brho}
c_1 {|\xi|}\leq\gau{K}(\xi)\leq c_2 {|\xi|}\,,\quad
\forall\xi\in\R^n.
\end{equation}

We say that $K\in\convopen$ is of class $C^2_+$
if $\partial K$ is of class $C^2$
and all the principal curvatures are strictly positive
functions on $\partial K$.
In this case, we define the $i$-th principal radius of
curvature at $x\in\partial K$ as
the reciprocal of the $i$-th principal curvature of $\partial K$
at $x$.
We remark that
if $K$ is of class $C^2_+$, then $K^0$ is also of
class $C^2_+$
(see \cite[p.~111]{Sch}).

Throughout the paper we shall assume that
\begin{equation}\label{f:ipoK}
K\in\convopen\
\textrm{is of class}\ C^2_+\,.
\end{equation}
Since $K$ will be kept fixed, from now on we shall use
the notation $\gauge = \gauge_K$ and $\pgauge=\pgauge_K$.

We collect here some known properties
of $\gauge$ and $\pgauge$
that will be frequently used
in the sequel.

\begin{thm}\label{t:sch}
Let $K$ satisfy $(\ref{f:ipoK})$. Then the following hold:
\par\noindent (i)
The functions $\gauge$ and $\pgauge$ are
convex, positively $1$-homoge\-neous
in $\R^n$, and
of class $C^2$ in $\R^n\setminus\{0\}$.
As a consequence,
\begin{gather*}
\gauge(t\, \xi) = t\, \gauge(\xi),\quad
D\gauge(t\, \xi) = D \gauge(\xi),\quad
D^2\gauge(t\, \xi) = \frac{1}{t}\, D^2\gauge(\xi),\\
\pgauge(t\, \xi) = t\, \pgauge(\xi),\
D\pgauge(t\, \xi) = D \pgauge(\xi),\
D^2\pgauge(t\, \xi) = \frac{1}{t}\, D^2\pgauge(\xi),
\end{gather*}
for every $\xi\in\R^n\setminus\{0\}$ and $t>0$.
Moreover
\begin{gather}
\pscal{D\gauge(\xi)}{\xi} = \gauge(\xi),\quad
\pscal{D\pgauge(\xi)}{\xi} = \pgauge(\xi)\,,
\quad\forall\xi\in\R^n\setminus\{0\}
\label{f:euler}\\
\label{f:eigenzero}
D^2\gauge(\xi)\, \xi = 0\,,\quad
D^2\pgauge(\xi)\, \xi = 0\,,
\qquad\forall \xi\in\R^n\setminus\{0\}\,.
\end{gather}
\par\noindent (ii)
For every $\xi,\eta\in\R^n$, we have
\begin{equation}\label{f:sublin}
\gauge(\xi+\eta)\leq \gauge(\xi) + \gauge(\eta),\quad
\pgauge(\xi+\eta)\leq \pgauge(\xi) + \pgauge(\eta),
\end{equation}
and equality holds if and only if $\xi$ and $\eta$ belong
to the same ray, that is,
$\xi = \lambda\, \eta$ or $\eta = \lambda\, \xi$
for some $\lambda\geq 0$.
\par\noindent (iii)
The ei\-gen\-val\-ues of the second differential
$D^2\gauge$ at $\uv\in S^{n-1}$
are $0$ (with corresponding eigenvector $\uv$)
and the principal radii of curvature of
$\partial K^0$ at the unique point $p\in\partial K^0$
at which $\uv$ is attained as an outward normal vector.
Symmetrically, the ei\-gen\-val\-ues of
$D^2\pgauge$ at $\uv\in S^{n-1}$
are $0$ (with corresponding eigenvector $\uv$)
and the principal radii of curvature of
$\partial K$ at the unique point $p\in\partial K$
at which $\uv$ is attained as an outward normal vector.
\end{thm}

\section{Distance from the boundary}

Throughout the paper, we shall assume that
\begin{equation}\label{f:Omega}
\Omega\subset\R^n\
\textrm{ is a nonempty, bounded, open connected set
of class $C^2$}.
\end{equation}

Let us define the function
\begin{equation}\label{f:d}
\dist(x) = \inf_{y\in\partial\Omega} \pgauge(x-y),
\qquad x\in\overline{\Omega},
\end{equation}
that measures the distance from
the boundary $\partial\Omega$ to a point $x\in\overline{\Omega}$
in the Minkowski norm
associated to the polar function $\pgauge$ of $\gauge$.
Since $\partial\Omega$ is a compact subset of $\R^n$
and $\pgauge$ is a continuous function,
the infimum in the definition of $\dist$ is achieved.
We shall denote by $\proj(x)$ the set of projections of $x$
onto $\partial\Omega$, that is
\begin{equation}\label{f:Pi}
\proj(x) = \{y\in\partial\Omega;\ \dist(x) = \pgauge(x-y)\},
\qquad x\in\overline{\Omega}.
\end{equation}
By abuse of notation,
when $\proj(x) = \{x_0\}$ then we shall use $\proj(x)$
to indicate the point $x_0$.

It is well-known that $\dist$ is a viscosity solution
of the Hamilton-Jacobi equation
\begin{equation}\label{f:HJ}
\gauge(Du) = 1\qquad
\textrm{in}\ \Omega.
\end{equation}
More precisely, it is the unique viscosity solution
of (\ref{f:HJ}) satisfying the boundary condition
$u=0$ on $\partial\Omega$.
Moreover, $\dist(x) > 0$ for every $x\in\Omega$,
$\dist\in\Lip(\overline{\Omega})$,
and $\gauge(Du(x)) = 1$ for a.e.~$x\in\Omega$
(see \cite{BaCD,Li}).

\medskip
We say that $x\in\Omega$ is a
\textsl{regular point} of $\Omega$
if $\proj(x)$ is a singleton.
We say that $x\in\Omega$ is a
\textsl{singular point} of $\Omega$
if $x$ is not a regular point.
We denote by $\Sigma\subseteq\Omega$ the set
of all singular points of $\Omega$.
It is well-known that $\dist$ is differentiable
at every regular point of $\Omega$
(see \cite{BaCD,CaSi,Li}; see also
Theorem~\ref{t:cm6}(i) below).

{}From now on,
for every $x_0\in\partial\Omega$ we shall denote
by $\curv_1(x_0),\ldots, \curv_{n-1}(x_0)$ the principal
curvatures of $\partial\Omega$ at $x_0$,
and by $\nor(x_0)$ the \textsl{inward} normal unit vector
to $\partial\Omega$ at $x_0$.
We extend these functions to
$\overline{\Omega}\setminus {\Sigma}$
by setting
\[
\nor(x) = \nor(\proj(x)),\quad
\curv_i(x) = \curv_i(\proj(x)), \quad
i=1,\ldots,n-1,\
x\in\overline{\Omega}\setminus {\Sigma}\,.
\]

We collect in the following theorem
all the results proved in \cite{CMf} that are
relevant to the subsequent analysis.

\begin{thm}\label{t:cm6}
Let $\Omega$ and $K$ satisfy respectively $(\ref{f:Omega})$
and $(\ref{f:ipoK})$. Then the following hold.

\noindent
(i) $\overline{\Sigma}\subset\Omega$,
and the Lebesgue measure of
$\overline{\Sigma}$ is zero.

\noindent
(ii) 
Let $x\in\Omega$ and $x_0\in\proj(x)$.
Then, for every $z\in [x_0,x)$, $\dist$ is differentiable at $z$
and
\begin{equation}\label{f:diffd}
D\dist(z) = \frac{\nor(x_0)}{\gauge(\nor(x_0))}\,.
\end{equation}

\noindent
(iii)
The function $\dist$ is of class $C^2$ on
$\overline{\Omega}\setminus\overline{\Sigma}$.
\end{thm}

\begin{proof}
See  Remark~4.16, Corollary~6.9,
Lemma~4.3
and Theorem~6.10 in \cite{CMf}.
\end{proof}

At any point $x_0\in\partial\Omega$
there is a unique inward ``normal'' direction $p(x_0)$
with the properties
$\proj(x_0+t p(x_0)) = \{x_0\}$
and $\dist(x_0+t p(x_0)) = t$
for $t\geq 0$ small enough
(see \cite[Remark~4.5]{CMf}).
More precisely, these properties
hold true for $t\in [0, \tau(x_0))$,
where $\tau(x_0)$ is the
normal distance to the cut locus
$\overline{\Sigma}$, defined below
(see \cite[Proposition~4.8]{CMf}).
It can be proved that $p(x_0) = D\gauge(\nor(x_0))$
(see \cite[Lemma~2.2]{LN}
and \cite[Proposition~4.4]{CMf}).
{}From Theorem~\ref{t:cm6}(ii) and the positive $0$-homogeneity
of $D\gauge$ it is clear that
$D\gauge(\nor(x_0)) = D\gauge(D\dist(x_0))$.
Summarizing,
given $x\in\overline{\Omega}$ we have that
\begin{equation}\label{f:xpx}
x_0\in\proj(x)\quad\Longleftrightarrow\quad
x=x_0+t\, D\gauge(D\dist(x_0))\
\textrm{for some}\ t\in [0,\tau(x_0)],
\end{equation}
and, in such case, $\dist(x) = t$.

The above considerations motivate the following definition.

\begin{defn}
The \textsl{normal distance to cut locus}
of a point $x\in\overline{\Omega}$ is
defined by
\begin{equation}\label{f:tau}
\tau(x) =
\begin{cases}
\min\{t\geq 0;\
x + t D\gauge(D\dist(x))\in\overline{\Sigma}\},
&\textrm{if $x\in\overline{\Omega}\setminus\overline{\Sigma}$},\\
0,
&\textrm{if $x\in\overline{\Sigma}$}.
\end{cases}
\end{equation}
The \textsl{cut point} $\cut(x)$
of $x\in\overline{\Omega}\setminus\overline{\Sigma}$ is defined by
$\cut(x) = x + \tau(x)D\gauge(D\dist(x))$.
\end{defn}

\begin{prop}\label{l:taupos}
Let $\Omega$ satisfy $(\ref{f:Omega})$.
Then $\tau$ is continuous in $\overline{\Omega}$.
Furthermore,
there exists $\mu > 0$ such that
$\tau(x_0) \geq \mu$ for every $x_0\in\partial\Omega$.
\end{prop}

\begin{proof}
See \cite{CMf}, Lemma~4.1 and Theorem~6.7.
\end{proof}


{}From Theorem~\ref{t:cm6}(iii), the function $\dist$
is of class $C^2$ on
$\overline{\Omega}\setminus\overline{\Sigma}$.
We can then define the function
\begin{equation}\label{f:W}
W(x) = -D^2\gauge(D\dist(x))\, D^2\dist(x)\,,\qquad
x\in\overline{\Omega}\setminus\overline{\Sigma}\,.
\end{equation}
For any $x_0\in\partial\Omega$, let $T_{x_0}$ denote
the tangent space to $\partial\Omega$ at $x_0$.
If $x\in\overline{\Omega}\setminus{\Sigma}$
and $\proj(x) = \{x_0\}$,
we set $T_x = T_{x_0}$.
Observe that, by (\ref{f:diffd}) and (\ref{f:eigenzero}),
we have $D^2\gauge(D\dist(x))\, \nor(x) = 0$.
Then, for every $v\in T_x$,
we have $W(x)\, v\in T_x$.
Hence, we can define the map
\begin{equation}\label{f:bw}
\bw(x)\colon T_x\to T_x,\quad
\bw(x)\, w = W(x)\, w,
\end{equation}
that can be identified with a linear application from
$\R^{n-1}$ to $\R^{n-1}$.

We shall use the following results
(see \cite{CMf}, Lemmas~4.10 and~5.1).

\begin{lem}\label{l:W1}
Let $x_0\in\partial\Omega$.
Then
\[
\det [I-t\, W(x_0)] = \det[I_{n-1}-t\,\bw(x_0)]
\]
for every $t\in\R$, and both determinants are strictly
positive for every $t\in [0,\tau(x_0))$.
\end{lem}

\begin{rem}
Let $x_0\in\partial\Omega$.
We recall that
$D\gauge(\nor(x_0))$ is an eigenvector of
$W(x_0)$ with corresponding eigenvalue zero
(see \cite[Lemma~4.18]{CMf}).
On the other hand, from Lemma~\ref{l:W1}
we deduce that a number $\curv\neq 0$ is an
eigenvalue of $W(x_0)$ if and only if
it is an eigenvalue of $\bw(x_0)$.
\end{rem}

Although
the matrix $\bw(x_0)$ is not in general symmetric,
its eigenvalues
are real numbers,
and so its eigenvectors are real
(see \cite{CMf}, Remark~5.3).
The eigenvalues of $\bw(x_0)$
have an important geometric interpretation.

\begin{defn}[$\gauge$-curvatures]\label{d:curv}
Let $x_0\in\partial\Omega$.
The \textsl{principal $\gauge$-curvatures} of $\partial\Omega$ at $x_0$,
with respect to the Minkowski norm $\dist$,
are the eigenvalues
$\curvg_1(x_0)\leq\cdots\leq\curvg_{n-1}(x_0)$
of $\bw(x_0)$.
The corresponding eigenvectors are the
\textsl{principal $\gauge$-directions} of $\partial\Omega$
at $x_0$.
\end{defn}

\medskip
Up to now we have analyzed some properties of the matrices
$W(x_0)$ and $\bw(x_0)$ at points $x_0\in\partial\Omega$.
Now we are interested in the evolution of these matrices
along the transport ray starting from $x_0$.

\begin{lem}\label{l:W3}
Let $x_0\in\partial\Omega$, and define
\begin{equation}\label{f:defV}
V(t) = W(x_0+t\, D\gauge(D\dist(x_0)))\,,\quad
\bv(t) = \bw(x_0+t\, D\gauge(D\dist(x_0)))
\end{equation}
for $t\in [0,\tau(x_0))$. Then
\begin{equation}\label{f:V}
V(t)\, [I-t\, V(0)] = V(0),\qquad
\bv(t)\, [I_{n-1}-t\bv(0)] = \bv(0),
\end{equation}
for every $t\in [0,\tau(x_0))$.
Furthermore, $\trace V(t) = \trace \bv(t)$
for every $t\in [0,\tau(x_0))$.
\end{lem}

\begin{proof}
Let us consider the principal coordinate system at $x_0$,
i.e.\ the coordinate system such that $x_0 = 0$, $e_n = \nor(x_0)$
and $e_i$ coincides with the $i$-th principal direction
of $\partial\Omega$ at $x_0$, $i=1,\ldots,n-1$.
Let $X\colon\U\to\R^n$ be a local parametrization of
$\partial\Omega$ in a neighborhood of $x_0=0$.
The relation (\ref{f:diffd}) can be written an
\begin{equation}\label{f:dcar}
D\dist(X(y)+t\, D\gauge(D\dist(X(y)))) = \frac{N(y)}{\gauge(N(y))}\,,
\ y\in\U,\
t\in [0,\tau(X(y))),
\end{equation}
where $N(y) = \nor(X(y))$.
Recall that 
\begin{equation}\label{f:N}
N(0) = e_n,\quad
\frac{\partial N}{\partial y_i}(0) = -\curv_i e_i,\quad
i=1,\ldots,n-1,
\end{equation}
where $\curv_1,\ldots,\curv_{n-1}$ are the principal curvatures
of $\partial\Omega$ at $x_0$.
Differentiating (\ref{f:dcar}) with respect to $y_i$
at $y=0$ and using (\ref{f:N})
we obtain,
for every $i=1,\ldots,n-1$,
\begin{equation}\label{f:dd1}
D^2\dist(x_0+t D\gauge(D\dist(x_0)))
\left[I + t\, D^2\gauge(D\dist(x_0)) D^2\dist(x_0))
\right] e_i = -\xi_i,
\end{equation}
where
$\xi_i = \curv_i\left[\gauge(\nor)\, e_i-\pscal{D\gauge(\nor)}{e_i}\,\nor
\right]/\gauge(\nor)^2$,
and $\nor = \nor(x_0)$. 
Differentiating (\ref{f:dcar}) with respect to $t$ at $y=0$ we get
\begin{equation}\label{f:dd2}
D^2\dist(x_0+t\, D\gauge(D\dist(x_0)))\,
D\gauge(D\dist(x_0)) = 0\,.
\end{equation}
Let us apply
$D^2\gauge(D\dist(x_0+t\, D\gauge(D\dist(x_0)))) = D^2\gauge(D\dist(x_0))$
to both sides of equations (\ref{f:dd1}) and (\ref{f:dd2}).
Recalling the definition (\ref{f:W}) of $W$,
we obtain the relations
\[
\begin{split}
& W(x_0+t\, D\gauge(D\dist(x_0)))\,
\left[I-t\, W(x_0)\right]\, e_i = \overline{\xi}_i,\quad
i=1,\ldots,n-1,\\
& W(x_0+t\, D\gauge(D\dist(x_0)))\, D\gauge(D\dist(x_0)) = 0\,,
\end{split}
\]
where $\overline{\xi}_i = D^2\gauge(D\dist(x_0))\, \xi_i$.
We have that, for every $t\in [0,\tau(x_0))$,
\begin{equation}\label{f:const4}
\begin{cases}
V(t)\, [I-t\, V(0)]\, e_i = \overline{\xi}_i,\quad
i=1,\ldots,n-1,\\
V(t)\, [I-t\, V(0)]\,D\gauge(D\dist(x_0)) = 0\,,
\end{cases}
\end{equation}
where in the second identity we have used the fact that
$V(0)\,D\gauge(D\dist(x_0)) = 0$
(see 
\cite[Lemma~4.18]{CMf}).
Since $\pscal{D\gauge(D\dist(x_0))}{\nor} = \gauge(\nor) > 0$ from
the positive $1$-homogeneity of $\gauge$,
it follows that
the vectors $e_1,\ldots,e_{n-1}, D\gauge(D\dist(x_0))$
span $\R^n$, hence from (\ref{f:const4}) we infer that
the matrix $V(t)\, [I-t\, V(0)]$ is independent of $t$.
Evaluating this matrix at $t=0$, we finally obtain
the first identity in (\ref{f:V}).
{}From (\ref{f:diffd}) and (\ref{f:eigenzero}) we have that
\[
\pscal{V(t)\, v}{\nor} =
-\gauge(\nor) \pscal{D^2\dist(x_0+t D\gauge(D\dist(x_0)))\, v}%
{D^2\gauge(\nor)\,\nor} = 0
\]
for every $v\in\R^n$,
hence the second
identity in (\ref{f:V}) is also satisfied,
and $\trace V(t) = \trace\bv(t)$.
\end{proof}

\begin{rem}
Let $\bv(t)$ be the function defined in (\ref{f:defV}).
By definition, the eigenvalues
$\curvg_1,\ldots,\curvg_{n-1}$
of $\bv(0)$ are
the principal $\gauge$-curvatures of $\partial\Omega$
at $x_0$,
and the corresponding eigenvectors
$w_1,\ldots,w_{n-1}$ are the principal $\gauge$-directions
of $\partial\Omega$ at $x_0$.
{}From the identity (\ref{f:V}) we obtain that
\[
\bv(t)\, (1-t\, \curvg_i) w_i = \curvg_i\, w_i,\qquad
i=1,\ldots,n-1,\
t\in [0,\tau(x_0)).
\]
Since for every $t\in [0,\tau(x_0))$ the point
$x_0+t\, D\gauge(D\dist(x_0))$ belongs to
$\overline{\Omega}\setminus\overline{\Sigma}$,
then $1-t\,\curvg_i>0$
(see \cite[Lemma~5.4]{CMf}),
and hence the eigenvalues of $\bv(t)$ are
\[
\curvg_i(t) = \frac{\curvg_i}{1-t\,\curvg_i}\,,
\qquad
i=1,\ldots,n-1,
\]
with corresponding eigenvectors $w_1,\ldots,w_{n-1}$.
\end{rem}

\begin{prop}\label{p:M}
For every $x_0\in\partial\Omega$ let us define the function
\begin{equation}\label{f:defM}
M_{x_0}(s,t) =
\exp\left(-\int_s^t \trace\bw(x_0+\sigma\, D\gauge(D\dist(x_0))\, d\sigma\right),
\end{equation}
for $s,t\in [0,\tau(x_0))$.
Then
\begin{equation}\label{f:M}
\begin{split}
M_{x_0}(s,t) & = \frac{\det\bw(x_0+s\,D\gauge(\nor(x_0)))}%
{\det\bw(x_0+t\,D\gauge(\nor(x_0)))} =
\frac{\det\left[I_{n-1}-t\bw(x_0)\right]}
{\det\left[I_{n-1}-s\bw(x_0)\right]}
\\ &
=\prod_{i=1}^{n-1} \frac{1-t\,\curvg_i(x_0)}{1-s\,\curvg_i(x_0)}
\end{split}
\end{equation}
for every $s,t\in [0,\tau(x_0))$,
where $\curvg_1(x_0),\ldots,\curvg_{n-1}(x_0)$
are the principal $\gauge$-curvatures of $\partial\Omega$
at $x_0$.
\end{prop}

\begin{proof}
Let $\bv(t)$, $t\in [0,\tau(x_0))$,
be the function defined in (\ref{f:defV}).
{}From Lemmas~\ref{l:W1} and~\ref{l:W3} we have
\begin{equation}\label{f:V2}
\bv(t) = \bv(0)\,[I_{n-1}-t\bv(0)]^{-1},
\qquad
\forall t\in [0,\tau(x_0))\,.
\end{equation}
This implies that
the matrix-valued function $\bv(t)$ satisfies the
differential equation
\[
\bv'(t) = \bv(t)^2,
\qquad
\forall t\in [0,\tau(x_0))\,,
\]
and hence the function
$\delta(t) = \det\bv(t)$
is a solution of the differential equation
$\delta'(t) = \trace[\bv(t)]\, \delta(t)$
in $[0,\tau(x_0))$.
Now, the first equality in (\ref{f:M}) follows from the
fact that
\[
\delta(t) = \delta(s)\,
\exp\left(\int_s^t \trace\bv(\sigma)\, d\sigma\right)
= \delta(s)/ M_{x_0}(s,t)\,,
\]
for every $s,t\in [0,\tau(x_0))$.
The second equality follows from~(\ref{f:V2}).
The last equality is a direct consequence of the
fact that
$\curvg_1(x_0),\ldots,\curvg_{n-1}(x_0)$
are the eigenvalues of the matrix
$\bv(0)$. 
\end{proof}

\begin{prop}\label{p:M2}
The function $M_{x_0}(s,t)$, defined in
(\ref{f:defM}), is jointly continuous with respect to
$x_0\in\partial\Omega$ and $s,t\in [0,\tau(x_0))$.
Furthermore,
\begin{equation}\label{f:Mx3}
0\leq M_{x_0}(s,t) \leq
\prod_{i=1}^{n-1} (1+ T\, \widetilde{K}_{-})\,,
\qquad
\forall x_0\in\partial\Omega,\
0\leq s\leq t < \tau(x_0),
\end{equation}
where
\begin{gather}
T=\max\{\tau(x);\ x\in\partial\Omega\},\label{f:T}\\
\widetilde{K}_{-} = \max\{[\curvg_i(x)]_{-};\
x\in\partial\Omega,\ i=1,\ldots,n-1\},
\label{f:K}
\end{gather}
being
$[a]_{-} = \max\{0, -a\}$ the negative part
of a real number $a$.
\end{prop}

\begin{proof}
The continuity of $M_{x_0}(s,t)$ follows from
its definition and the continuity of $\curvg_i$.
The estimate (\ref{f:Mx3}) follows from
the representation formula (\ref{f:M}) and the estimate
\[
\frac{1-t\,\curvg_i(x_0)}{1-s\,\curvg_i(x_0)}
\leq 1+(t-s){[\curvg_i(x_0)]}_{-}
\]
that holds for every
$x_0\in\partial\Omega$
and $0\leq s\leq t < \tau(x_0)$.
\end{proof}

For every $x\in\overline\Omega\setminus\overline{\Sigma}$
let us define the function
\begin{equation}\label{f:defMx}
M_{x}(t) =
\exp\left(-\int_0^t \trace\bw(x+\sigma\, D\gauge(D\dist(x))\, d\sigma\right),\quad
t\in [0,\tau(x))\,,
\end{equation}
where $\bw$ is the matrix defined in (\ref{f:bw}).
For an explicit computation of $M_x$ it can be of some aid
to recall that
\[
\trace\bw(x) = \trace W(x) =
-\trace\left[D^2\gauge(D\dist(x))\, D^2\dist(x)\right]
\]
for every $x\in\overline\Omega\setminus\overline{\Sigma}$
(see Lemma~\ref{l:W3}).
Given $x\in\overline\Omega\setminus\overline{\Sigma}$,
let $\proj(x)=\{x_0\}$.
By (\ref{f:xpx})
we have that
\[
\bw(x+\sigma\, D\gauge(D\dist(x)) =
\bw(x_0+(\dist(x)+\sigma)\, D\gauge(D\dist(x)),
\qquad
\sigma\in [0,\tau(x)),
\]
which implies the relation
\begin{equation}\label{f:Mx}
M_{x}(t) =
M_{x_0}(\dist(x), \dist(x)+t),\qquad
\forall t\in [0,\tau(x))\,.
\end{equation}
{}From the identity (\ref{f:M}) we have that
\begin{equation}\label{f:Mx2}
M_{x}(t) =
\prod_{i=1}^{n-1} \frac{1-(\dist(x)+t)\,\curvg_i(x)}{1-\dist(x)\,\curvg_i(x)}
\end{equation}
where $\curvg_i(x) := \curvg_i(x_0)$, $i=1,\ldots,n-1$,
are the $\gauge$-curvatures of $\partial\Omega$ at $x_0$.

\section{Existence of solutions}
\label{s:exi}

In this section we recall the existence result
for system (\ref{f:syst1})-(\ref{f:syst2})
proved in \cite[Thm.~7.2]{CMf}.
The rigorous meaning of solution is the following.

\begin{defn}\label{d:sol}
Let $\Omega\subset\R^n$ satisfy (\ref{f:Omega})
and let $f\in \Cb(\Omega)$ be a
nonnegative function.
A solution of system (\ref{f:syst1})-(\ref{f:syst2})
is a pair $(u,v)$ of functions
satisfying the following properties:
\begin{itemize}
\item[1.]
$u\in\Lip(\overline{\Omega})$, $v\in \Cb(\Omega)$,
$u, v\geq 0$ in $\Omega$;
\item[2.]
$u = 0$ on $\partial\Omega$,
$\gauge(Du)\leq 1$ a.e.~in $\Omega$,
and $u$ is a viscosity solution of
\[
\gauge(Du) = 1\qquad
\textrm{in}\ \{v>0\}\,;
\]
\item[3.]
$v$ is a solution in the sense of distributions
of $-\dive(v\, D\gauge(Du)) = f$ in $\Omega$,
that is
\begin{equation}\label{f:weak}
\int_{\Omega} v(x)\,\pscal{D\gauge(Du(x))}{D\varphi(x)}\, dx
= \int_{\Omega} f(x)\, \varphi(x)\, dx
\end{equation}
for every $\varphi\in C^{\infty}_c(\Omega)$.
\end{itemize}
\end{defn}

\begin{rem}\label{r:dens}
Since $v\in L^{\infty}(\Omega)$
and $D\gauge(Du)\in [L^{\infty}(\Omega)]^n$,
by a standard density argument
(\ref{f:weak}) holds for every test function
$\varphi$ 
in the Sobolev space $W^{1,1}_0(\Omega)$.
\end{rem}

For $f\in \Cb(\Omega)$, let us define
the function
\begin{equation}\label{f:vf}
\vf(x) =
\begin{cases}
\displaystyle\int_0^{\tau(x)} f(x+t\, D\gauge(D\dist(x)))\, M_x(t)\, dt\,,
&\textrm{if $x\in\overline\Omega\setminus\overline{\Sigma}$},\\
0,
&\textrm{if $x\in\overline{\Sigma}$}\,.
\end{cases}
\end{equation}

\begin{thm}[Existence]\label{t:existence}
Let $\Omega\subset\R^n$ satisfy $(\ref{f:Omega})$
and let $f\geq 0$ be a bounded continuous function in $\Omega$.
Then, the pair $(\dist, \vf)$ is a solution to
$(\ref{f:syst1})$-$(\ref{f:syst2})$
in the sense of Definition~\ref{d:sol}.
\end{thm}

A rigorous proof of
Theorem~\ref{t:existence}
was given in \cite[Theorem~7.2]{CMf}.
In Section~\ref{s:uni} we shall prove that the pair
$(\dist, \vf)$ is essentially the unique
solution to problem (\ref{f:syst1})-(\ref{f:syst2}).
In order to gain some insight
in the representation formula (\ref{f:vf}),
a formal derivation of (\ref{f:vf}) might be in order.

Assume that $(\dist,v)$ is a solution of
(\ref{f:syst1})-(\ref{f:syst2}),
and that $v\in C^1(\Omega\setminus\overline{\Sigma})$,
with $v$ vanishing on $\overline{\Sigma}$.
Outside $\overline{\Sigma}$, the equation
$-\dive(v\, D\gauge(Du)) = f$ is satisfied pointwise, that is
\[
v(x)\, \trace W(x)-\pscal{Dv(x)}{D\gauge(D\dist(x))} = f(x)\,,
\]
where $W(x) = -D^2\gauge(D\dist(x))\, D^2\dist(x)$.
Furthermore, from Lemma~\ref{l:W3}
we have that $\trace W(x) = \trace\bw(x)$,
where $\bw(x)$ is the matrix defined in (\ref{f:bw}).

Let $x\in\Omega\setminus\overline{\Sigma}$, and
define
$\bar{v}(t) = v(x+t\, D\gauge(D\dist(x)))$,
$t\in [0,\tau(x)]$.
The function $\bar{v}(t)$ satisfies the following
linear differential equation
\[
\bar{v}'(t) = \left[\trace\bw(x+t\,D\gauge(D\dist(x)))\right]\,
\bar{v}(t) - f(x+t\,D\gauge(D\dist(x)))
\]
in $[0,\tau(x)]$, supplemented by the boundary condition
\[
\bar{v}(\tau(x)) = 0
\]
since $x+\tau(x)\,D\gauge(D\dist(x))\in\overline{\Sigma}$.
The solution of this Cauchy problem,
evaluated at $t=0$, gives
\[
v(x) = \bar{v}(0) =
\int_0^{\tau(x)} f(x+t\,D\gauge(D\dist(x)))\,
M_x(t)
\, dt\,,
\]
that is, the solution $v(x)$ has to be the function
defined in formula (\ref{f:vf}).

With this heuristic in mind,
our aim will be to prove that,
if $(u,v)$ is a solution to
(\ref{f:syst1})-(\ref{f:syst2}),
then $(\dist, v)$ is a solution too
(see Lemma~\ref{l:unique1}(ii)),
and that if $(\dist, v)$ is a solution
to (\ref{f:syst1})-(\ref{f:syst2}),
then $v$ must vanish on $\overline{\Sigma}$
(see Proposition~\ref{vzero}).
The first goal will be achieved using the
same arguments proposed in \cite{CCCG},
whereas the second one needs a new regularity
result for solutions of elliptic equations,
which seems to be of some interest by itself,
and that will be proved in the
following section.

\section{A regularity result}
\label{s:reg}

The aim of this section is to prove the following regularity
result.

\begin{thm}\label{t:reg}
Assume that $\gauge$ is the gauge function of a
convex body $K$ satisfying (\ref{f:ipoK}).
Let $A\subset\R^n$ be an open bounded set,
and let $u\in W^{1,\infty}(A)$
be a solution in the sense of distributions
of the equation
\begin{equation}\label{f:equa2}
-\dive(D\gauge(Du)) = 0\qquad \textrm{in $A$},
\end{equation}
that is
\[
\int_A \pscal{D\gauge(Du(x))}{D\varphi(x)}\, dx=0
\]
for every $\varphi\in C^{\infty}_c(A)$. If in addition
\begin{equation}\label{f:equa1}
\gauge(Du(x)) = 1\qquad \textrm{a.e.\ in $A$,}
\end{equation}
then $u\in C^{1,\alpha}_{loc}(A)$.
\end{thm}

We recall the standard regularity result about
solutions to the equation
\begin{equation}\label{f:dgeq}
-\dive(a(Du))=0 \qquad \textrm{in $A$},
\end{equation}
(see \cite[\S4.6]{LU}, \cite[\S8.2]{Gi}).
\begin{thm}
Assume that the vector--valued function $a$ belongs to $C^1(A)$
and satisfies the following growth conditions: there exist $p>1$,
and $\alpha_0,\beta_0>0$ such that for every $\xi\in \R^n$
\begin{align}
& |a(\xi)|+(1+|\xi|^2)^{\frac{1}{2}}|Da(\xi)| \leq
\alpha_0 (1+|\xi|^2)^{\frac{p-1}{2}}\,, \label{f:standard1} \\
& \pscal {Da(\xi) w}{w} \geq \beta_0 (1+|\xi|^2)^{\frac{p-2}{2}}|w|^2\,,\ \forall w\in\R^n\,.
\label{f:standard3}
\end{align}
Then every solution $u\in W^{1,\infty}(A)$ of (\ref{f:dgeq}) belongs to $C^{1,\alpha}_{loc}(A)$.
\end{thm}

In our case, $\gauge$ is a positively $1$-homogeneous function
of class $C^2$ in $\R^n\setminus\{0\}$ satisfying the bounds
(\ref{f:brho}), and then $D\gauge$ is a positively $0$-homogeneous
function of class $C^1$ in $\R^n\setminus\{0\}$, but in general $\gauge$ is not
even differentiable
at the origin. Moreover, (\ref{f:standard3}) cannot be verified
near the origin. Hence we have no chance to apply directly the
standard regularity results to solutions of equation (\ref{f:equa2}).

An easy trick in order to have the ``right" growth is to
consider the function
\begin{equation}\label{f:gamma}
\gamma(\xi) = \frac{1}{2}\, \gauge(\xi)^2,\qquad \xi\in\R^n\,.
\end{equation}
We have that
\begin{equation}\label{f:dgamma}
D\gamma(\xi) =
\begin{cases}
\gauge(\xi)\, D\gauge(\xi),
&\textrm{if $\xi\neq 0$,}\\
0,
&\textrm{if $\xi=0$.}
\end{cases}
\end{equation}

Clearly, a function $u\in W^{1,\infty}(A)$ satisfying
(\ref{f:equa1})
is a solution to
(\ref{f:equa2})
if and only if it is a solution to
\begin{equation}\label{f:sistb}
\begin{cases}
-\dive(D\gamma(Du)) = 0
&\textrm{in $A$},\\
\gauge(Du) = 1
&\textrm{a.e.\ in $A$}.
\end{cases}
\end{equation}
Moreover, $\gamma(\xi)$ is a positively $2$--homogeneous function
of class $C^2$ in $\R^n\setminus \{0\}$, hence $D^2\gamma(\xi)$
is 0--homogeneous and continuous in $\R^n\setminus \{0\}$.
In particular
the matrix--valued function $D^2\gamma$ is bounded in $\R^n\setminus\{0\}$, and
\begin{equation}\label{f:estd2ga}
\|D^2\gamma\|\leq c_5 \doteq\max_{i,j=1,\ldots,n}
\max\{|D^2_{i,j}\gamma(\xi)|,\ \xi\in S^{n-1}\}\,.
\end{equation}

The following positive constants will be used
throughout this section:
\begin{equation}\label{f:const}
\begin{split}
c_3 & = \max\{|D\gauge(\xi)|;\ \xi\in \R^n\setminus\{0\}\}
=\max\{|D\gauge(\uv)|;\ \uv\in S^{n-1}\},\\
c_4 &
=\max\{|D^2\gauge(\uv)|;\ \uv\in S^{n-1}\},\\
r^0 & = \min\{r^0_i(p);\ p\in\partial K^0,\ i=1,\ldots,n-1\},\\
R^0 & = \max\{r^0_{i}(p);\ p\in\partial K^0,\ i=1,\ldots,n-1\}\,
\end{split}
\end{equation}
where $r^0_1(p)\leq \ldots\leq r^0_{n-1}(p)$ are the principal curvatures
of $\partial K^0$ at $p$.
We remark that $r^0>0$ since
$K^0$ is of class $C^2_+$.

The first technical tool is to prove that $\gamma$ satisfies some
growth conditions  similar to (\ref{f:standard1})
and (\ref{f:standard3}) with $p=2$.

\begin{lem}\label{l:reg}
Let $\gamma$ be the function defined in $(\ref{f:gamma})$.
Then
\begin{align}
& |D\gamma(\xi)| \leq c_2 c_3\,|\xi|\,,
\label{f:reg1}\\
& \pscal{D\gamma(\xi)}{\xi} \geq c_1^2 |\xi|^2\,,
\label{f:reg2}
\end{align}
for every $\xi\neq 0$. Here $c_1$, $c_2$, $c_3$ are the constants defined in
$(\ref{f:brho})$ and $(\ref{f:const})$. Moreover
there exists a constant $c_6>0$, independent of $\xi\neq 0$, such that
\begin{equation}\label{f:reg3}
\pscal{D^2\gamma(\xi)w}{w} \geq c_6 |w|^2,
\end{equation}
for every $w\in\R^n$.
\end{lem}

\begin{proof}
Since $\gamma$ is 2--homogeneous, by Euler's formula $\pscal{D\gamma(\xi)}{\xi}=
2\gamma(\xi)$, and then
\[
\pscal{D\gamma(\xi)}{\xi} =
\gauge(\xi)^2 \geq c_1^2 |\xi|^2\,.
\]
On the other hand, by (\ref{f:dgamma}),
\[
|D\gamma(\xi)| \leq \gauge(\xi) |D\gauge(\xi)| \leq c_2 c_3\,|\xi|\,.
\]
It remains to prove that there exists a constant
$c_6 > 0$ such that
\begin{equation}\label{f:stg}
\pscal{D^2\gamma(\xi)w}{w}\geq c_6 |w|^2,\quad
\forall w\in\R^n,\ \forall \xi\neq 0\,.
\end{equation}
Since $\gamma$ is a convex function, of class
$C^2$ in $\R^n\setminus\{0\}$, we have that
the quadratic form $w\mapsto \pscal{D^2\gamma(\xi)w}{w}$
is positive semidefinite for every $\xi\neq 0$.
We shall show that, in fact, it is positive definite
uniformly with respect to $\xi\neq 0$.

Fixed $\xi\neq 0$ and using the notation $\uv=\xi/|\xi|$, we have
\[
D^2\gamma(\xi)=D^2\gamma(\uv) = \gauge(\uv)\, D^2\gauge(\uv) +
D\gauge(\uv) \otimes D\gauge(\uv),
\]
so that
\[
\pscal{D^2\gamma(\xi)w}{w}=\pscal{D^2\gamma(\uv)w}{w}=
\gauge(\uv)\, \pscal{D^2\gauge(\uv)w}{w} +
\pscal{D\gauge(\uv)}{w}^2
\]
for every $w\in\R^n$.

Fixed $w\in\R^n$, let us denote by $\lambda=\pscal{\uv}{w}$, and by
$\wbot$ the projection of $w$ on the orthogonal space
$\mathcal{L}^\bot(\uv)$ to $\uv$,
so that $w=\wbot+\lambda \uv$, $\pscal{\wbot}{\uv}=0$,
and $|w|^2=|\wbot|^2+\lambda^2$.

{}From Theorem~\ref{t:sch} we have that
$D^2\gauge(\uv)\uv=0$ and $D^2\gauge(\uv)$ is positive definite in
$\mathcal{L}^\bot(\uv)$ with
\[
\pscal{D^2\gauge(\uv)\wbot}{\wbot}\geq r^0|\wbot|^2\,,
\]
for every $\wbot\in \mathcal{L}^\bot(\uv)$,
where $r^0>0$ is the constant defined in (\ref{f:const}).
Hence we obtain
\[
\pscal{D^2\gauge(\uv)w}{w}=
\pscal{D^2\gauge(\uv)\wbot}{\wbot}\geq r^0|\wbot|^2\,,
\]
for every $w\in\R^n$.
On the other hand, by the 1--homogeneity of $\gauge$ we have
\[
\pscal{D\gauge(\uv)}{\wbot+\lambda \uv}=
\pscal{D\gauge(\uv)}{\wbot}+\lambda\gauge(\uv)\,.
\]
Hence we get the inequality
\begin{align*}
\pscal{D^2\gamma(\xi)w}{w} & \geq  r^0\gauge(\uv)|\wbot|^2+
\pscal{D\gauge(\uv)}{\wbot}^2
+\lambda^2\gauge(\uv)^2+2\lambda\gauge(\uv)\pscal{D\gauge(\uv)}{\wbot}\,,
\end{align*}
for every $w=\wbot+\lambda \uv\in\R^n$, $\wbot\in \mathcal{L}^\bot(\uv)$.

It remains to prove that there exists $c_6>0$ independent of $\uv$ such that
\begin{align*}
 r^0\gauge(\uv)|\wbot|^2 & +\pscal{D\gauge(\uv)}{\wbot}^2 +
\lambda^2\gauge(\uv)^2
+2\lambda\gauge(\uv)\pscal{D\gauge(\uv)}{\wbot}
\geq  c_6\left(|\wbot|^2+\lambda^2\right)\,.
\end{align*}
We have that
\begin{equation}\label{f:magg7}
\begin{split}
\frac{r^0}{2}\gauge(\uv)|\wbot|^2+\pscal{D\gauge(\uv)}{\wbot}^2 & \geq
\left(1+ \frac{r^0}{2}\frac{\gauge(\uv)}{|D\gauge(\uv)|^2}\right)
\pscal{D\gauge(\uv)}{\wbot}^2
\\ & \geq
\left(1+\frac{r^0c_1}{2c_3^2}\right)\pscal{D\gauge(\uv)}{\wbot}^2
\end{split}
\end{equation}
where $c_1>0$ and $c_3>0$ are defined respectively
in (\ref{f:brho}) and (\ref{f:const}).
If $c=1+\frac{r^0c_1}{2c_3^2}$, from (\ref{f:magg7}) we get
\begin{align*}
 r^0\gauge(\uv) & |\wbot|^2  +\pscal{D\gauge(\uv)}{\wbot}^2 +
\lambda^2\gauge(\uv)^2
+2\lambda\gauge(\uv)\pscal{D\gauge(\uv)}{\wbot} \\
 & \geq \frac{r^0}{2}\gauge(\uv)|\wbot|^2 + c
 \left(\pscal{D\gauge(\uv)}{\wbot}+\frac{\lambda\gauge(\uv)}{c}\right)^2+
 \left(1-\frac{1}{c}\right)\lambda^2\gauge(\uv)^2 \\
 & \geq\frac{r^0c_1}{2}\,|\wbot|^2+ \left(1-\frac{1}{c}\right)c_1^2\lambda^2
 \geq c_6\left(|\wbot|^2 +\lambda^2\right)\,,
 \end{align*}
 where $c_6=\min\left(\frac{r^0c_1}{2},\left(1-\frac{1}{c}\right)c_1^2\right)>0$.
\end{proof}

Thanks to the estimates (\ref{f:reg1}), (\ref{f:reg2}),
and (\ref{f:reg3})
we can now prove that a solution $u$ to the first
equation in (\ref{f:sistb}) belongs to $H^2_{loc}(A)$.
This part of the proof is based on
a standard argument in regularity theory (see e.g.\
\cite[\S8.2]{Gi}).
The only point that should be stressed concerns the
regularity of $\gamma$. Namely, $\gamma$ is of class
$C^2$ in $\R^n\setminus\{0\}$, but in general
the positively $0$-homogeneous function $D^2\gamma$
is not even defined at the origin.
In our case this is not a real problem, since the
condition $\gauge(Du(x)) = 1$ guarantees that $Du(x)$ stays
always outside a ball centered at the origin.

\begin{lem}\label{l:akka2}
Every solution $u\in W^{1,\infty}(A)$ of (\ref{f:sistb}) belongs to
$H^2_{loc}(A)$.
\end{lem}

\begin{proof}
The function $u$ solves
\begin{equation}\label{f:eqdeb}
\int_A \pscal {D\gamma(D u(x))}{D\varphi(x)}\, dx =0
\end{equation}
for every $\varphi \in C^{\infty}_c(A)$.
Let us fix $\varphi$,
a coordinate direction $e_i$, and $h\in\R$ such that
$\varphi_h(x)=\varphi(x-he_i)$ has compact support in $A$.
Choosing $\varphi_h$ as a test function in (\ref{f:eqdeb})
and making a change of variables in the integral we obtain
\begin{equation}\label{f:eqdeb1}
\int_A \pscal {D\gamma(D u(x+he_i))}{D\varphi(x)}\, dx =0\,.
\end{equation}
Taking the difference between (\ref{f:eqdeb}) and (\ref{f:eqdeb1}),
we get
\begin{equation}\label{f:eqdeb2}
\int_A \pscal {D\gamma(D u(x+he_i))-D\gamma(Du(x))}{D\varphi(x)}\, dx =0\,.
\end{equation}
On the other hand we have that $\gauge(Du(x)) =\gauge(Du(x+he_i))= 1$
for a.e.~$x\in A$.
Then the function $\alpha(t)=D\gamma((1-t)Du(x)+t Du(x+h e_i))$ is continuous in
the interval $[0,1]$.
In addition, either $\alpha\in C^1([0,1])$,
or there exists $t_0\in(0,1)$ such that
$\alpha$ is of class $C^1$ and with bounded derivatives in $[0,1]\setminus \{t_0\}$.
Hence $\alpha(t)$ is a
Lipschitz function, and
\begin{equation}\label{f:fu}
\begin{split}
& D \gamma(D u(x+he))-D\gamma(Du(x)) \\
& =
\int_0^1 \frac{d}{dt}[D\gamma((1-t)Du(x)+t Du(x+h e_i))]\, dt \\
&=\int_0^1 D^2\gamma((1-t)Du(x)+t Du(x+h e_i))[Du(x+h e_i))-Du(x))]dt,
\end{split}
\end{equation}
for a.e.\ $x\in A$. Finally, the matrix
\[
L_h(x)= \int_0^1 D^2\gamma((1-t)Du(x)+t Du(x+h e_i))\, dt
\]
is well defined and, for $w\in\R^n$, the integrals of the kind
\[
\pscal{L_h(x)w}{w} = \int_0^1 \pscal{D^2\gamma((1-t)Du(x)+t Du(x+h e_i)) w}{w}\, dt
\]
satisfy for a.e.\ $x\in A$
\begin{align}
& c_6|w|^2\leq\pscal{L_h(x)w}{w}\,, \quad \forall w\in\R^n,\label{f:Lbounds1}\\
& \|L_h(x)\|\leq c_5\,, \label{f:Lbounds2}
\end{align}
where $c_6$ and $c_5$ are the positive constants defined in Lemma~\ref{l:reg}
in (\ref{f:estd2ga}).

If we denote by
$\Delta_hu(x)=\frac{u(x+he_i)-u(x)}{h}$,
by (\ref{f:fu}) equation (\ref{f:eqdeb2})
can be rewritten as
\begin{equation}\label{f:eqdeb3}
\int_{A}\pscal{L_h(x) D \Delta_hu(x)}{D\varphi(x)}\, dx=0
\end{equation}
for every $\varphi\in C^{\infty}_c(A)$ and for every
$|h|<\textrm{dist}(\textrm{supp}\varphi,\partial A)$.
By a density
argument, we have that (\ref{f:eqdeb3}) remains valid for every test function in
$H^1_0(A')$,
where $A'$ is an open set compactly contained in $A$, and for every
$|h|<\textrm{dist}(A',\partial A)$.
Hence we can choose $\varphi=\eta^2\Delta_hu $ as test function in (\ref{f:eqdeb3}),
where
$\eta\in C^{\infty}_c(A)$ is defined in the following way: given $x_0\in A$
and $r>0$
such that the ball $B_{2r}(x_0)$ is compactly contained in $A$, we require that
$0\leq\eta(x)\leq 1$ in $A$, $\eta(x)=0$ in $A \setminus B_{2r}(x_0)$, $\eta(x)=1$
in $B_{r}(x_0)$, and there exists $m>0$ such that $|D\eta(x)|\leq m/r$ in $A$.
With this choice of the test function, (\ref{f:eqdeb3}) becomes
\begin{equation}\label{f:eqdeb6}
\begin{split}
\int_{A}& \pscal{L_h(x) D\Delta_hu(x)}{D\Delta_h u(x)}\,\eta^2(x)\, dx
\\ &
=
-2\int_{A}\pscal{L_h(x) D\Delta_hu(x)}{D\eta(x)}\,\eta(x)\Delta_hu(x)\, dx\,,
\end{split}
\end{equation}
for every $|h|<\textrm{dist}(x_0,\partial A)-2r$.
Recalling (\ref{f:Lbounds1}) we have that
\begin{equation}\label{f:eqdeb4}
c_6\int_A|D\Delta_hu(x)|^2\eta^2\, dx\leq
\int_{A} \pscal{L_h(x) D\Delta_hu(x)}{D\Delta_h u(x)}\eta^2(x)\, dx\,.
\end{equation}
On the other hand,
by (\ref{f:Lbounds2}) and Young's inequality we obtain
that there exists a constant $C>0$ such that
\begin{equation}\label{f:eqdeb5}
\begin{split}
& \left|\int_{A}  \pscal{L_h(x) D\Delta_hu(x)}{D\eta(x)}\, \eta(x)\Delta_hu(x)\, dx
\right|
\\
& \leq C \left(\varepsilon \int_A|D\Delta_hu(x)|^2\eta^2(x)\, dx + \frac 1\varepsilon
\int_A|D\eta(x)|^2|\Delta_hu(x)|^2\, dx\right)\,,
\end{split}
\end{equation}
for every $\varepsilon>0$. Choosing $\varepsilon$ small enough, from
(\ref{f:eqdeb6}), (\ref{f:eqdeb4}), (\ref{f:eqdeb5}), and the estimates
of $|D\eta|$, we get
\begin{equation}
\int_{B_{r}(x_0)}|D\Delta_hu(x)|^2\, dx\leq
\frac M{r^2}\int_{B_{2r}(x_0)}|\Delta_hu(x)|^2\, dx\,,
\end{equation}
which implies, by a standard argument, that $u\in H_{loc}^2(A)$
(see Lemmas~7.23 and~7.24
in \cite{GT}, or \cite[\S8.1]{Gi} ).
\end{proof}

\begin{proof}[Proof of Theorem \ref{t:reg}]
We already know that a solution $u\in W^{1,\infty}(A)$ of
(\ref{f:equa1})--(\ref{f:equa2}) is also a solution to
(\ref{f:sistb}).
Fixed $\varphi\in C^{\infty}_c(\Omega)$, we can choose
$\frac{\partial \varphi}{\partial x_i}$ as test
function in the weak formulation (\ref{f:eqdeb}), obtaining
\begin{equation}
\int_A \pscal {D\gamma(Du(x))}{D\frac{\partial \varphi}{\partial x_i}(x)}\, dx
=0
\end{equation}
Moreover, by Lemma \ref{l:akka2}, the function $u$ belongs to $H_{loc}^2(A)$.
Then an integration by parts leads
\begin{equation}
\int_A \pscal {D^2\gamma(Du(x))D\frac{\partial u}{\partial x_i}(x)}
{D\varphi(x)}\, dx=0\,.
\end{equation}
for every $\varphi\in C^{\infty}_c(\Omega)$.
Then the partial derivative $\frac{\partial u}{\partial x_i}$
is a bounded solution of the linear elliptic equation
\begin{equation}
-\dive(D^2\gamma(Du(x))Dv)= 0 \qquad \textrm{in \ }A\,,
\end{equation}
where the matrix $D^2\gamma(Du(x))$ satisfies the hypothesis
of the De Giorgi--Nash regularity result (see \cite[\S3.14]{LU}, Theorem 14.1,
\cite{GT}, Theorem 8.22).
Hence we can conclude that the partial derivatives
of $u$ are locally H\"older continuous.
\end{proof}

\begin{rem}
Let us define the function $\gamma_p(\xi) = \gauge(\xi)^p/p$, $p>1$.
Clearly a function $u\in W^{1,\infty}(A)$ is a solution to
(\ref{f:equa1})-(\ref{f:equa2})
if and only if it is a solution to
\begin{equation}\label{f:sistc}
\begin{cases}
-\dive(D\gamma_p(Du(x))) = 0
&\textrm{in $A$},\\
\gauge(Du(x)) = 1
&\textrm{in $A$}
\end{cases}
\end{equation}
where, as usual, the first equation is interpreted in the sense
of distributions and the second in viscosity sense.
Arguing as in the proof of Lemma~\ref{l:reg},
it can be proved that
there exists a positive constant $c_p$ such that, for every $\xi\neq 0$,
\begin{align*}
& |D\gamma_p(\xi)| \leq c_2^{p-1} c_3\,|\xi|^{p-1}\,,\\
& \pscal{D\gamma_p(\xi)}{\xi} \geq c_1^p |\xi|^p\,,\\
& \pscal{D^2\gamma_p(\xi)w}{w} \geq c_p |\xi|^{p-2} |w|^2,\quad\forall w\in\R^n.
\end{align*}
We remark that, if $p>2$, then $\gamma_p$ is of class $C^2$ on $\R^n$,
and the estimates above hold for every $\xi\in\R^n$. On the other hand, if $p>2$, we cannot
obtain an estimate of the type (\ref{f:standard3})
near the origin, due to the $p$--homogeneity of the function $\gamma_p$,
which implies the $(p-2)$--homogeneity of $D^2\gamma_p(\xi)$.
For this reason we have considered the case $p=2$.
\end{rem}

\section{Uniqueness}
\label{s:uni}

This section is devoted to the proof of the following
uniqueness result.

\begin{thm}\label{t:unique}
Let $(u,v)$ be a solution of system
$(\ref{f:syst1})$-$(\ref{f:syst2})$ in the sense
of Definition~\ref{d:sol}.
Then $v=\vf$, where $\vf$ is the function
defined in $(\ref{f:vf})$,
and $u = \dist$
in $\Omega_f = \{x\in\Omega;\ \vf(x) > 0\}$.
\end{thm}

The proof of Theorem \ref{t:unique}
is essentially based on the
techniques developed in
\cite{Egan,CaCa,Pr,CCCG}.
We will first prove the uniqueness of the first component of the solution of system
(\ref{f:syst1})--\eqref{f:syst2}.
More precisely, we will show that if $(\uz,\vz)$ is a solution of system
(\ref{f:syst1})--\eqref{f:syst2},
then $\uz = \dist$ in $\Omega_f:=\{x\in\Omega~:~\vf(x)>0\}$
(see Proposition~\ref{unique2} below).

Let us consider the 
functional
$\Phi\colon H_0^1(\Omega)\times L^2_+(\Omega)\to\R$,
where $L^2_+(\Omega) = \{\vv\in L^2(\Omega); \ \vv\geq 0\}$,
defined by
\begin{equation}\label{phi}
\Phi(\uu,\vv)=-\int_\Omega f(x)\uu(x)\; dx+
\int_\Omega \vv(x)\, \left[
\gauge(D\uu(x))-1\right]\, dx.
\end{equation}

\begin{lem}\label{saddle-point}
If $(\uz,\vz)$ is a solution of system \eqref{f:syst1}--\eqref{f:syst2},
then $(\uz,\vz)$ is a saddle point of
$\Phi$, in the sense that
$$
\Phi(\uz,\vv)\leq\Phi(\uz,\vz)\leq \Phi(\uu,\vz)\qquad \forall
(\uu,\vv)\in
H_0^1(\Omega)\times
L^2_+(\Omega).
$$
\end{lem}

\begin{proof}
Since $(\uz,\vz)$ is a solution of (\ref{f:syst1})--\eqref{f:syst2}, then
$$
\int_\Omega \vz(x)\, \left[
\gauge(D\uz(x))-1\right]\, dx=0
$$
and
$$
\int_\Omega \vv(x)
\left[
\gauge(D\uz(x))-1\right]\, dx\leq 0, \quad \forall
\vv\in L^2_+(\Omega).
$$
Hence, for any $\vv\in L^2_+(\Omega)$ we have
\begin{equation}\label{saddle1}
\begin{split}
\Phi(\uz,\vz) & =
-\int_\Omega f(x)\uz(x)\, dx
\\ & \geq
-\int_\Omega
f(x)\uz(x)\, dx +
\int_\Omega \vv(x)\, \left[
\gauge(D\uz(x))-1\right]\, dx
= \Phi(\uz,\vv).
\end{split}
\end{equation}
Moreover, by the convexity of $\gauge$,
for any $\uu\in H^1_0(\Omega)$ we have
\begin{equation}\label{f:triag}
\gauge(D\uu(x)) - \gauge(D\uz(x))\geq
\pscal{D\gauge(D\uz(x))}{D\uu(x)-D\uz(x)}
\end{equation}
for a.e.~$x\in\Omega$.
By Remark~\ref{r:dens},
we can choose $\varphi = \uu-\uz\in H^1_0(\Omega)$
as test function in (\ref{f:weak}), obtaining
\[
 -\int_\Omega f(x)(\uu(x)-\uz(x))\, dx
+\int_\Omega \vz(x)\pscal{D\gauge(D\uz(x))}{D\uu(x)-D\uz(x)}\, dx=0.
\]
Thus, by (\ref{f:triag}), for any $\uu\in H^1_0(\Omega)$,
\begin{equation}\label{saddle2}
\begin{split}
\Phi(\uu,\vz)  -\Phi(\uz,\vz)
\geq&
-\int_\Omega f(x)(\uu(x)-\uz(x))\, dx
\\ & +\int_\Omega \vz(x)\pscal{D\gauge(D\uz(x))}{D\uu(x)-D\uz(x)}\, dx
= 0\,.
\end{split}
\end{equation}
Collecting together (\ref{saddle1})
and (\ref{saddle2}) we get the conclusion.
\end{proof}

In what follows we shall use the
set of functions
\begin{equation}\label{f:lip}
\Lip^1_{\gauge}(\Omega):=
\{\uu\in \Lip(\overline{\Omega});\
\uu = 0\ \textrm{on}\ \partial\Omega,\
\gauge(D\uu)\leq 1\ \textrm{a.e.~in}\ \Omega\}.
\end{equation}
It can be checked that $\uu\in\Lip^1_{\gauge}(\Omega)$
if and only if
$\uu\in \Lip(\overline{\Omega})$, $\uu = 0$ on $\partial\Omega$,
and
\begin{equation}\label{f:lip1}
\uu(x) - \uu(y) \leq \pgauge(x-y)
\qquad
\forall x,y\in\overline{\Omega},\
\textrm{with}\ [x,y]\subset\overline{\Omega}
\end{equation}
(see \cite[Chap.~6]{Li}).

\begin{lem}\label{l:unique1}
If $(u,v)$ is a solution of  system
$(\ref{f:syst1})$--$(\ref{f:syst2})$, then
the following hold.
\begin{itemize}
\item[(i)]
$u = \dist$ in $\spt(f)$.
\item[(ii)]
$(\dist,v)$ is a
solution of $(\ref{f:syst1})$--$(\ref{f:syst2})$.
\end{itemize}
 \end{lem}

\begin{proof}
(i) By the maximality property of viscosity solutions
we have that
$\uz\leq\dist$ in $\Omega$.
On the other hand,
by Lemma~\ref{saddle-point},
$\Phi(\uz,\vz)\leq\Phi(\uu,\vz)$
for any $\uu\in \Lip^1_{\gauge}(\Omega)$.
In addition, we have
\[
\begin{split}
\Phi(\uz,\vz) & = -\int_\Omega f(x)\uz(x)\, dx,\\
\Phi(\uu,\vz) & = -\int_\Omega f(x)\uu(x)\, dx
+ \int_\Omega \vz(x)\left[\gauge(D\uu(x))-1\right]\, dx
\\ & \leq -\int_\Omega f(x)\uu(x)\, dx\,.
\end{split}
\]
Then
\[
\int_\Omega f(x)\uu(x)\, dx \leq
\int_\Omega f(x)\uz(x)\, dx
\qquad
\forall \uu\in \Lip^1_{\gauge}(\Omega)\,.
\]
Choosing $\uu = \dist$, we obtain that
$\uz = \dist$ on $\spt(f)$.

(ii) From~(i)
and Lemma~\ref{saddle-point}
we have that
\[
\Phi(\dist,v)=
\Phi(u,v)\leq \Phi(\uu,v).
\qquad
\forall \uu\in H^1_0(\Omega).
\]
Hence, for every test function
$\varphi\in C^{\infty}_c(\Omega)$ and every $h>0$ we have that
\[
\begin{split}
0 & \leq \Phi(\dist+h\,\varphi, v)-\Phi(\dist,v)
\\ & =
-h\int_\Omega f(x)\varphi(x)\, dx +
\int_\Omega v(x)\left[\gauge(D(\dist(x)+h\, \varphi(x))-\gauge(D\dist(x))
\right]\, dx\,.
\end{split}
\]
Since $\gauge$ is convex
we have
\[
\left|\frac{\gauge(D(\dist(x)+h\, \varphi(x)))-\gauge(D\dist(x))}{h}
\right|
\leq \|D\gauge\|_{\infty}\,
\|D\varphi\|_{L^{\infty}(\Omega)}
\]
for a.e.~$x\in\Omega$ and every $h> 0$, where
\[
\|D\gauge\|_{\infty} :=
\sup_{\xi\neq 0} |D\gauge(\xi)|
= \max_{\xi\in S^{n-1}} |D\gauge(\xi)|
\]
due to the positive $0$-homogeneity of $D\gauge$.
Hence, from the differentiability of $\gauge$ in $\R^n\setminus\{0\}$
and the dominated convergence theorem we get
\[
0\leq -\int_\Omega f(x)\varphi(x)\, dx
+\int_\Omega v(x)\, \pscal{D\gauge(D\dist(x))}{D\varphi(x)}\, dx\,.
\]
Replacing $\varphi$ by $-\varphi$ we also get the opposite inequality.
\end{proof}

\begin{prop}\label{unique2}
If $(u,v)$ is a solution of  system $(\ref{f:syst1})$--$(\ref{f:syst2})$,
then $u = \dist$ in the set
$\Omega_f = \{ x\in \Omega;\ \vf(x)>0\}$,
where $\vf$ is the function defined by $(\ref{f:vf})$.
\end{prop}

\begin{proof}
Let $x\in\Omega_f\subseteq\Omega\setminus\overline{\Sigma}$.
By the definition (\ref{f:vf}) of $\vf$,
and taking into account that
$M_x(t)>0$ for every $t\in [0,\tau(x))$
and $f\geq 0$,
we deduce that there exists $t_0\in (0,\tau(x))$
such that, at the point
$x_0 = x + t_0\, D\gauge(D\dist(x))$,
one has $f(x_0) > 0$.
Since
\[
\dist(x_0) = \dist(x) + t_0 = \dist(x) + \pgauge(x_0-x)\,,
\]
from this identity, Lemma~\ref{l:unique1}(i),
(\ref{f:lip1}) and
the inequality $u\leq\dist$ we get
\[
\dist(x) = \dist(x_0) - \pgauge(x_0-x)
= u(x_0) - \pgauge(x_0-x)
\leq u(x) \leq\dist(x)\,,
\]
hence $u(x) = \dist(x)$.
\end{proof}

Now that the uniqueness of the first component of the solution of system
(\ref{f:syst1})--\eqref{f:syst2} is proved, it remains to prove the uniqueness of the second one.
In order to do so, we
will first exhibit for such a function a representation formula on the set
$\Omega\setminus\overline{\Sigma}$ and then analyze its behavior on
$\overline{\Sigma}$.

\begin{prop}\label{vrepr}
If $(\dist,v)$ is a solution of  system $(\ref{f:syst1})$--$(\ref{f:syst2})$,
then for any $z_0\in \Omega\setminus \overline{\Sigma}$
and $\theta\in (0,\tau(z_0))$ we have
\begin{equation}\label{repr}
\begin{split}
v(z_0) & - v(z_0+\theta\, D\gauge(D\dist(z_0)))\, M_{z_0}(\theta)
\\ &
= \int_0^\theta f(z_0+t\, D\gauge(D\dist(z_0)))\, M_{z_0}(t)
\, dt\,.
\end{split}
\end{equation}
\end{prop}

\begin{proof}
Let $z_1 = z_0+\theta\, D\gauge(D\dist(z_0))$,
let $\proj(z_0)=\{x_0\}$,
and define $t_0 = \dist(z_0) = \pgauge(z_0-x_0)$,
$t_1 = \dist(z_1) = \pgauge(z_1-x_0) = t_0+\theta$.
Let $Y\colon\U\to\R^n$, $\U\subset\R^{n-1}$ open,
a local parametrization of $\partial\Omega$ in a
neighborhood of $x_0$, such that $Y(0) = x_0$.
Let $\Psi\colon\U\times\R\to\R^n$ be the map
\begin{equation}\label{f:psi}
\Psi(y,t) = Y(y) + t D\gauge(\nor(Y(y))),
\quad (y,t)\in \U\times\R\,.
\end{equation}
Choose $r>0$ such that
$U(r) := \{y\in\R^{n-1};\ |y|\leq r\}\subset\U$, and
\[
D(r) := \{\Psi(y,t);\ y\in U(r),\
t\in [t_0,t_1]\} \subset \Omega\setminus\overline{\Sigma}\,.
\]
The set $D(r)$ can be viewed as a tubular neighborhood
of the segment $[z_0, z_1]$.
Let us define
\[
S_i(r) = \{\Psi(y, t_i);\ y\in U(r)\}\,,
\qquad i=0,1,
\]
and let $S_2(r)$ denote the lateral surface of $D(r)$, i.e.\
\[
S_2(r) = \{\Psi(y,t);\ y\in\partial U(r),\
t\in [t_0,t_1]\}\,.
\]
All these surfaces are of class $C^1$ and
are oriented with the outward normal with
respect to $D(r)$.

For $\epsilon>0$ small enough let
$\psi_{\epsilon}\colon\R\to\R$ and
$\eta_{\epsilon}\colon\R^{n-1}\to\R$
be the functions defined by
\[
\psi_{\epsilon}(t) =
\begin{cases}
0 &
\textrm{if $t\leq t_0$ or $t\geq t_1$,}\\
1 &
\textrm{if $t\in [t_0+\epsilon, t_1-\epsilon]$,}\\
\frac{t-t_0}{\epsilon} &
\textrm{if $t\in (t_0, t_0+\epsilon)$,}\\
\frac{t_1-t}{\epsilon} &
\textrm{if $t\in (t_1-\epsilon, t_1)$.}
\end{cases}\quad
\eta_{\epsilon}(y) =
\begin{cases}
0 & \textrm{if $|y|\geq r$,}\\
1 & \textrm{if $|y|\leq r-\epsilon$,}\\
\frac{r-|y|}{\epsilon} &
\textrm{if $r-\epsilon < |y|< r$.}
\end{cases}
\]
Let $\varphi_{\epsilon}$ be the function defined by
\[
\varphi_{\epsilon}(x):=
\begin{cases}
\psi_{\epsilon}(t) \eta_{\epsilon}(y),
&\textrm{if $\exists\ y\in U(r)$ and $t\in [t_0,t_1]$ s.t.\ $x=\Psi(y,t)$},\\
0,
&\textrm{otherwise}.
\end{cases}
\]
It is clear that $\varphi_{\epsilon}$
belongs to $\Lip(\overline{\Omega})$ and has
support contained in $D(r)$, hence can be used as test function
in (\ref{f:weak}).

It is plain that $\varphi_{\epsilon}$ converges monotonically to
$1$ in the interior of $D(r)$ as $\epsilon\to 0^+$, hence
\begin{equation}\label{f:eqxf}
\lim_{\epsilon\to 0^+}\int_{\Omega} f\, \varphi_{\epsilon}\, dx
= \int_{D(r)} f\, dx\,.
\end{equation}
Let us compute the right-hand side of (\ref{f:weak})
when $\varphi = \varphi_{\epsilon}$
and $u=\dist$.
On $D(r)$ the test function $\varphi_{\epsilon}$ is
defined by the relation
\[
\varphi_{\epsilon}(\Psi(y,t)) = \psi_{\epsilon}(t)\,
\eta_{\epsilon}(y),\qquad
y\in U(r),\ t\in [t_0, t_1]\,.
\]
Differentiating the relation above with respect to $t$
and recalling the definition (\ref{f:psi}) of $\Psi$
we obtain
\[
\pscal{D\varphi_{\epsilon}(\Psi(y,t))}%
{D\gauge(D\dist(\Psi(y,t)))} = \psi'_{\epsilon}(t)\,
\eta_{\epsilon}(y),\qquad
y\in U(r),\ t\in [t_0, t_1]\,.
\]
Then, taking into account that $\psi'_{\epsilon}(t) = 0$
for $t\in (t_0+\epsilon, t_1-\epsilon)$, we get
\begin{equation}\label{f:eqxx}
\begin{split}
& \int_{\Omega}  v\, \pscal{D\gauge(D\dist)}{D\varphi_{\epsilon}}\, dx
\\ & =
\int_{U(r)}\int_{t_0}^{t_1} v(\Psi)\,
\pscal{D\gauge(D\dist(\Psi))}{D\varphi_{\epsilon}(\Psi)}
\, \det D\Psi\, dt\, dy
\\ & =
\int_{U(r)}\int_{t_0}^{t_1} v(\Psi)\,
\psi'_{\epsilon}(t)\,
\eta_{\epsilon}(y)
\, \det D\Psi\, dt\, dy
\\ &
=
I_0(\epsilon)+I_1(\epsilon)+I_2(\epsilon)+
I_3(\epsilon)\,,
\end{split}
\end{equation}
where
\[
\begin{split}
I_0(\epsilon) & = \frac{1}{\epsilon}\int_{U(r-\epsilon)}
\int_{t_0}^{t_0+\epsilon}
v(\Psi)\, \det D\Psi\, dt\, dy\,,\\
I_1(\epsilon) & = -\frac{1}{\epsilon}\int_{U(r-\epsilon)}
\int_{t_1-\epsilon}^{t_1}
v(\Psi)\, \det D\Psi\, dt\, dy\,,\\
I_2(\epsilon) &= \frac{1}{\epsilon}\int_{U(r)\setminus U(r-\epsilon)}
\int_{t_0}^{t_0+\epsilon}
v(\Psi)\, \eta_{\epsilon}(y)\,\det D\Psi\, dt\, dy\,,\\
I_3(\epsilon) &= -\frac{1}{\epsilon}\int_{U(r)\setminus U(r-\epsilon)}
\int_{t_1-\epsilon}^{t_1}
v(\Psi)\, \eta_{\epsilon}(y)\,\det D\Psi\, dt\, dy\,.
\end{split}
\]
Since $v$ and $D\Psi$ are bounded in $D(r)$,
an explicit computation leads to
$|I_2(\epsilon)+I_3(\epsilon)| \leq C \epsilon$.
Passing to the limit in (\ref{f:eqxx}),
by the continuity of $v(\Psi)$ and $D\Psi$
we obtain
\[
\lim_{\epsilon\to 0^+}
\int_{\Omega} v\, \pscal{D\dist}{D\varphi_{\epsilon}}\, dx =
\sum_{i=0}^1 (-1)^i
\int_{U(r)} v(\Psi(y,t_i))\, \det D\Psi(y,t_i)\, dy\,.
\]
Recalling (\ref{f:eqxf}),
we finally obtain
\begin{equation}\label{iden}
\int_{D(r)} f(x)\, dx =
\sum_{i=0}^1 (-1)^i\int_{U(r)} v(\Psi(y,t_i))\,
\det D\Psi(y, t_i) \,dy\,.
\end{equation}

As a last step we want to pass to the limit
as $r\to 0^+$.
{}From the continuity of $v$ and $D\Psi$ we get
\begin{equation}\label{est2}
\begin{split}
\lim_{r\to 0^+} & \frac{1}{\sigma(r)}\int_{U(r)}
v(\Psi(y,t_i))\,
\det D\Psi(y, t_i) \,dy
=
v(z_i)\,
\det D\Psi(0, t_i)\,,
\end{split}
\end{equation}
where $\sigma(r)=\omega_{n-1}r^{n-1}$
is the area of the ball with radius $r>0$ in $\R^{n-1}$.
Finally
\begin{equation}\label{est6}
\begin{split}
\lim_{r\to 0} & \frac{1}{\sigma(r)}\int_{D(r)} f(x)\,dx
\\ &
=
\lim_{r\to 0}\frac{1}{\sigma(r)}
\int_{U(r)} \int_{t_0}^{t_1} f(\Psi(y,t))\, \det D\Psi(y,t)\, dt\, dy
\\ & =
\int_{t_0}^{t_1} f(\Psi(0,t))\,
\det D\Psi(0, t)\, dt
\\ & =
\int_{0}^{\theta} f(z_0+t\,D\gauge(D\dist(z_0)))\,
\det D\Psi(0, t_0+t)\, dt\,.
\end{split}
\end{equation}
{}From Lemma~4.10 in \cite{CMf}
we have that
\[
\det D\Psi(0,t) =
\sqrt{G}\,\gauge(\nor(x_0))\,\det(I_{n-1}-t\, \bw(x_0)),
\qquad
t\in [t_0, t_1]\,,
\]
where $G$ is the determinant of the matrix of the
metric coefficients.
Collecting together (\ref{est2}), (\ref{est6}),
recalling the identity (\ref{iden}),
and dividing by
$\det D\Psi(0,t_0)$
we obtain
\[
\begin{split}
& v(z_0)
-v(z_1)\,
\frac{\det(I_{n-1}-t_1\, \bw(x_0))}{\det(I_{n-1}-t_0\, \bw(x_0))}
\\ & =
\int_{0}^{\theta} f(z_0+t\,D\gauge(D\dist(z_0)))\,
\frac{\det(I_{n-1}-(t_0+t)\, \bw(x_0))}{\det(I_{n-1}-t_0\, \bw(x_0))}\, \, dt\,.
\end{split}
\]
The representation formula (\ref{repr})
now follows from (\ref{f:M}) and the definition
(\ref{f:Mx}) of $M_{z_0}$.
\end{proof}

For the proof of Proposition~\ref{vzero} below
we need two more technical ingredients.
The first one is the regularity result proved
in Theorem~\ref{t:reg}.
The second one is the following convergence lemma due
to H.~Brezis (see \cite[Theorem~1]{Br}).

\begin{lem}\label{l:brezis}
Let $\gamma\colon\R^n\to\R$ be a strictly convex
function, satisfying the linear growth condition
\[
\gamma(\xi)\geq c_0 |\xi| - b_0,
\qquad\forall\xi\in\R^n,
\]
for some positive constants $b_0$ and $c_0$.
Let $(u_k)_k\subset [L^1(\Omega)]^n$ be a sequence
of functions converging to $u\in [L^1(\Omega)]^n$
in the weak $L^1$ topology,
and assume that $\gamma(u)$, $\gamma(u_k)\in [L^1(\Omega)]^n$
for every $k\in\N$.
If $\lim_k \int_\Omega \gamma(u_k)\, dx = \int_\Omega\gamma(u)\, dx$,
then $(u_k)_k$ converges to $u$ in the strong $L^1$ topology.
\end{lem}

\begin{prop}\label{vzero}
If $(\dist,v)$ is a solution of  system $(\ref{f:syst1})$--$(\ref{f:syst2})$,
then $v(x) =  0$ for every
$x\in\overline{\Sigma}$.
\end{prop}

\begin{proof}
Since $v$ is a continuous function, it suffices to prove
that $v = 0$ on $\Sigma$.
Let us fix any $x_0\in\Sigma$ and choose
$\epsilon>0$ sufficiently small such that
$B_\epsilon(x_0)\subset \Omega$.
Then, for any $x\in B_1(0)$ set
\[
d_\epsilon(x):=
\frac{\dist(x_0+\epsilon x)-\dist(x_0)}{\epsilon},\
v_\epsilon(x):=v(x_0+\epsilon x),\
f_\epsilon(x):=f(x_0+\epsilon x).
\]
By construction, for any $\epsilon>0$ as above $d_\epsilon(0)=0$ and
\[
\gauge(Dd_\epsilon(x))=\gauge(D\dist(x_0+\epsilon x))=1
\qquad\textrm{for a.e.}\ x\in B_1(0)\,.
\]
Hence, there exist a sequence
$(\epsilon_j)_j$,
$\epsilon_j\to 0^+$ and a
Lipschitz function
$d_0\colon B_1(0)\to\R$
such that 
$(d_{\epsilon_j})_j$ converges to $d_0$ uniformly in $B_1(0)$.
Moreover, since
$\gauge(Dd_{\epsilon_j}(x))=1$ in the viscosity sense in $B_1(0)$,
by \cite[Proposition 2.2]{BaCD} also
$\gauge(Dd_0(x))=1$ in the
viscosity sense in $B_1(0)$,
which gives $\gauge(Dd_0(x))=1$ almost everywhere.
Since $(Dd_{\epsilon_j})_j$ is bounded in $L^{\infty}$,
we can also assume that it converges to $Dd_0$
in the weak $L^1$ topology.

\noindent
{}From (\ref{f:reg3}) in Lemma~\ref{l:reg},
we have that the function $\gamma(\xi) = \gauge(\xi)^2/2$
is strictly convex in $\R^n$.
Since $\gamma(Dd_{\epsilon_j}(x))=1/2$ and $\gamma(Dd_0(x)) = 1/2$
for a.e.~$x\in B_1(0)$, we have that
\[
\lim_{j\to\infty} \int_{B_1(0)} \gamma(Dd_{\epsilon_j}(x))\, dx=
\frac{\omega_n}{2}
=\int_{B_1(0)} \gamma(Dd_{0}(x))\,dx\,.
\]
Recalling that
$(Dd_{\epsilon_j})_j$ converges to $Dd_0$
in the weak $L^1$ topology,
from Lemma~\ref{l:brezis} we conclude that
$(Dd_{\epsilon_j})_j$ converges to $Dd_0$
in the strong $L^1$ topology.

\noindent
Finally, the functions
$v_{\epsilon_j}$ and $f_{\epsilon_j}$ defined above uniformly converge to $v(x_0)$
and $f(x_0)$ respectively and the
pair $(d_{\epsilon_j},v_{\epsilon_j})$ solves
\begin{equation}\label{f:eqep}
-\dive (v_{\epsilon_j} D\gauge(Dd_{\epsilon_j}))=
\epsilon_j\, f_{\epsilon_j}\qquad \mbox{in } B_1(0)
\end{equation}
in the sense of distributions,
due to the fact that $(\dist,v)$ solves (\ref{f:syst1})--\eqref{f:syst2}.
Upon observing that $D\gauge(Dd_{\epsilon_j})$
converges to $D\gauge(Dd_0)$ in the strong $L^1$ topology,
we can pass to the
limit as $j\to \infty$
in (\ref{f:eqep}),
obtaining that
$d_0$ is a weak solution of
\[
-\dive (v(x_0)\, D\gauge(Dd_0(x)))=0\qquad x\in B_1(0).
\]
Now, if $v(x_0)\neq 0$,
then $d_0$ must be a solution to
\[
\begin{cases}
-\dive (D\gauge(Dd_0))=0
& \textrm{in $B_1(0)$},\\
\gauge(D d_0) = 1
& \textrm{a.e.~in $B_1(0)$}.
\end{cases}
\]
{}From Theorem~\ref{t:reg} we have that
$d_0\in C^{1,\alpha}(B_1(0))$.
On the other hand, $d_0$ cannot be differentiable in $x=0$, because $d_0$ is
the `blow up' of the distance function around
a singular point $x_0$.
Hence $v(x_0)=0$ and the proof is complete.
\end{proof}

The last two propositions allow us to prove Theorem~\ref{t:unique} as a
simple corollary.
Indeed, we already know by
Proposition \ref{unique2} that
if $(u,v)$ is a solution of system (\ref{f:syst1})--\eqref{f:syst2},
then $u = \dist$ on
the set $\Omega_f=\{\vf>0\}$.
So it only remains to prove that $v = \vf$ in $\Omega$,
where $\vf$ is given by (\ref{f:vf}).
Proposition \ref{vzero} guarantees that $v = 0$ in
$\overline{\Sigma}$, while Proposition \ref{vrepr} implies
that for any $z_0\in \Omega\setminus\overline{\Sigma}$ and
$\theta \in (0,\tau(z_0))$
\[
\begin{split}
v(z_0) & - v(z_0+\theta\, D\gauge(D\dist(z_0)))\, M_{z_0}(\theta)
\\ & = \int_0^\theta f(z_0+t\, D\gauge(D\dist(z_0)))\, M_{z_0}(t)
\, dt\,.
\end{split}
\]
Hence, letting $\theta \to \tau(z_0)^-$ and using the continuity of $v$
we obtain that $v(z_0) = \vf(z_0)$.

\medskip
\par\noindent
\textbf{Acknowledgements.}
The authors wish to thank Nicola Fusco for
fruitful discussions about the regularity result
proved in Section~\ref{s:reg},
and an anonymous referee for the careful reading of
the manuscript.

%

\end{document}